\documentclass{amsart}
\usepackage{amsfonts}
\usepackage{amsmath}
\usepackage{amssymb}
\usepackage{amscd}

\input{tcilatex}
\theoremstyle{definition}
\theoremstyle{remark}
\numberwithin{equation}{section}

\begin{document}
\author{D. S. Lubinsky}
\address{School of Mathematics\\
Georgia Institute of Technology\\
Atlanta, GA 30332-0160\\
USA.\\
lubinsky@math.gatech.edu}
\title[Universality Limits]{A New Approach to Universality Limits at the
Edge of the Spectrum}
\date{January 4, 2007}
\maketitle

\begin{abstract}
We show how localization and smoothing techniques can be used to establish
universality at the edge of the spectrum for a fixed positive measure $\mu $
on $\left[ -1,1\right] $. Assume that $\mu $ is a regular measure, and is
absolutely continuous in some closed neighborhood $J$ of $1$. Assume that in 
$J$, $\mu ^{\prime }\left( x\right) =h\left( x\right) \left( 1-x\right)
^{\alpha }\left( 1+x\right) ^{\beta }$, where $h\left( 1\right) >0$ and $h$
is continuous at 1. Then universality at $1$ for $\mu $ follows from
universality at $1$ for the classical Jacobi weight $\left( 1-x\right)
^{\alpha }\left( 1+x\right) ^{\beta }$.
\end{abstract}

\dedicatory{Dedicated to the 60th birthday of Percy Deift}

\section{Results\protect\footnote{%
Research supported by NSF grant DMS0400446 and US-Israel BSF grant 2004353}}

Let $\mu $ be a finite positive Borel measure on $\left( -1,1\right) $. Then
we may define orthonormal polynomials%
\begin{equation*}
p_{n}\left( x\right) =\gamma _{n}x^{n}+...,\gamma _{n}>0,
\end{equation*}%
$n=0,1,2,...$ satisfying the orthonormality conditions%
\begin{equation*}
\int_{-1}^{1}p_{n}p_{m}d\mu =\delta _{mn}.
\end{equation*}%
These orthonormal polynomials satisfy a recurrence relation of the form%
\begin{equation*}
xp_{n}\left( x\right) =a_{n+1}p_{n+1}\left( x\right) +b_{n}p_{n}\left(
x\right) +a_{n}p_{n-1}\left( x\right) ,
\end{equation*}%
where%
\begin{equation*}
a_{n}=\frac{\gamma _{n-1}}{\gamma _{n}}>0\text{ and }b_{n}\in \mathbb{R}%
\text{, }n\geq 1,
\end{equation*}%
and we use the convention $p_{-1}=0$. Throughout $w=\frac{d\mu }{dx}$
denotes the absolutely continuous part of $\mu $. A classic result of E.A.
Rakhmanov \cite{Simon2005} asserts that if $w>0$ a.e. in $\left[ -1,1\right] 
$, then $\mu $ belongs to the Nevai-Blumenthal class $\mathcal{M}$, that is 
\begin{equation*}
\lim_{n\rightarrow \infty }a_{n}=\frac{1}{2}\text{ and }\lim_{n\rightarrow
\infty }b_{n}=0.
\end{equation*}%
A class of measures that contains $\mathcal{M}$ is the class of regular
measures \cite{StahlTotik1992}, defined by the condition 
\begin{equation*}
\lim_{n\rightarrow \infty }\gamma _{n}^{1/n}=\frac{1}{2}.
\end{equation*}

One of the key limits in random matrix theory, the so-called universality
limit \cite{Deift1999}, involves the reproducing kernel%
\begin{equation*}
K_{n}\left( x,y\right) =\sum_{k=0}^{n-1}p_{k}\left( x\right) p_{k}\left(
y\right)
\end{equation*}%
and its normalized cousin 
\begin{equation*}
\widetilde{K}_{n}\left( x,y\right) =w\left( x\right) ^{1/2}w\left( y\right)
^{1/2}K_{n}\left( x,y\right) .
\end{equation*}%
In \cite{Lubinsky2006}, we presented a new approach to this universality
limit, proving:\newline
\newline
\textbf{Theorem 1.1} \newline
\textit{Let }$\mu $\textit{\ be a finite positive Borel measure on }$\left(
-1,1\right) $\textit{\ that is regular. Let} $I$\textit{\ be a closed
subinterval of }$\left( -1,1\right) $\textit{\ in which }$\mu $\textit{\ is
absolutely continuous.}\newline
\textit{(a) Assume that }$w$\textit{\ is positive and continuous in }$I$%
\textit{. Let }$\Omega $ \textit{denote the modulus of continuity of }$w$%
\textit{\ in }$I$\textit{, so that for }$\delta >0,$ 
\begin{equation*}
\Omega \left( w;\delta \right) =\sup \left\{ \left\vert w\left( x\right)
-w\left( y\right) \right\vert :x,y\in I\text{ and }\left\vert x-y\right\vert
\leq \delta \right\} .
\end{equation*}%
\textit{Assume that }$w$\textit{\ satisfies the Dini condition }%
\begin{equation*}
\int_{0}^{1}\frac{\Omega \left( w;t\right) }{t}dt<\infty .
\end{equation*}%
\textit{Then if }$I^{\prime }$\textit{\ is a compact subinterval of }$I^{0}$%
\textit{, we have}%
\begin{equation*}
\lim_{n\rightarrow \infty }\frac{\widetilde{K}_{n}\left( x+\frac{a}{%
\widetilde{K}_{n}\left( x,x\right) },x+\frac{b}{\widetilde{K}_{n}\left(
x,x\right) }\right) }{\widetilde{K}_{n}\left( x,x\right) }=\frac{\sin \pi
\left( a-b\right) }{\pi \left( a-b\right) }.
\end{equation*}%
\textit{uniformly for }$x\in I^{\prime }$\textit{\ and }$a,b$\textit{\ in
compact subsets of the real line.}\newline
\textit{(b) Assume that }$w$\textit{\ is bounded above and below by positive
constants, and moreover, }$w$\textit{\ is Riemann integrable in }$I$\textit{%
. Then if }$p>0$\textit{\ and} $I^{\prime }$ \textit{is a closed subinterval
of }$I^{0},$%
\begin{equation*}
\lim_{n\rightarrow \infty }\int_{I^{\prime }}\left\vert \frac{\widetilde{K}%
_{n}\left( x+\frac{a}{\widetilde{K}_{n}\left( x,x\right) },x+\frac{b}{%
\widetilde{K}_{n}\left( x,x\right) }\right) }{\widetilde{K}_{n}\left(
x,x\right) }-\frac{\sin \pi \left( a-b\right) }{\pi \left( a-b\right) }%
\right\vert ^{p}dx=0,
\end{equation*}%
\textit{\ uniformly for }$a,b$\textit{\ in compact subsets of the real line. 
\newline
}

We also established $L_{1}$ analogues assuming less on $w$. However, we have
subsequently realized that the smoothness condition in (a) can be dropped,
and all we need is that $w$ is continuous on $I$. The technique of \cite%
{Lubinsky2006} involved a Taylor series expansion in $\alpha ,\beta $ of $%
K_{n}\left( x+\frac{\alpha }{n},x+\frac{\beta }{n}\right) $, a localization
technique, and a smoothing technique. In this paper, we show how
localization and smoothing can be applied at the edge $1$ of the spectrum.
As far as the author is aware, the most general result to date for Jacobi
type weights is due to Kuijlaars and Vanlessen \cite{KuijlaarsVanlessen2002}%
. Let $\mu $ be absolutely continuous, and $w$ have the form%
\begin{equation}
w\left( x\right) =h\left( x\right) w^{\left( a,\beta \right) }\left(
x\right) =h\left( x\right) \left( 1-x\right) ^{\alpha }\left( 1+x\right)
^{\beta },
\end{equation}%
where $h$ is positive and analytic in $\left[ -1,1\right] $. They showed
that uniformly for $a,b$ in bounded subsets of $\left( 0,\infty \right) ,$
as $n\rightarrow \infty ,$ 
\begin{equation}
\frac{1}{2n^{2}}\tilde{K}_{n}\left( 1-\frac{a}{2n^{2}},1-\frac{b}{2n^{2}}%
\right) =\mathbb{J}_{\alpha }\left( a,b\right) +O\left( \frac{a^{\alpha
/2}b^{\alpha /2}}{n}\right) .
\end{equation}%
Here 
\begin{equation*}
\mathbb{J}_{\alpha }\left( u,v\right) =\frac{J_{\alpha }\left( \sqrt{u}%
\right) \sqrt{v}J_{\alpha }^{\prime }\left( \sqrt{v}\right) -J_{\alpha
}\left( \sqrt{v}\right) \sqrt{u}J_{\alpha }^{\prime }\left( \sqrt{u}\right) 
}{2\left( u-v\right) }
\end{equation*}%
is the Bessel kernel of order $\alpha $, and $J_{\alpha }$ is the usual
Bessel function of the first kind and order $\alpha $. Our result is:\newline
\newline
\textbf{Theorem 1.2}\newline
\textit{Let }$\mu $\textit{\ be a finite positive Borel measure on }$\left(
-1,1\right) $\textit{\ that is regular. Assume that for some }$\rho >0$, $%
\mu $ \textit{is absolutely continuous in }$J=\left[ 1-\rho ,1\right] $, 
\textit{and in }$J$\textit{, its absolutely continuous component has the
form }$w=hw^{\left( \alpha ,\beta \right) }$, \textit{where }$\alpha ,\beta
>-1.$ \textit{Assume that }$h\left( 1\right) >0$\textit{\ and }$h$\textit{\
is} \textit{continuous at }$1$. \textit{Then uniformly for }$a,b$\textit{\
in compact subsets of }$\left( 0,\infty \right) $\textit{, we have }%
\begin{equation}
\lim_{n\rightarrow \infty }\frac{1}{2n^{2}}\tilde{K}_{n}\left( 1-\frac{a}{%
2n^{2}},1-\frac{b}{2n^{2}}\right) =\mathbb{J}_{\alpha }\left( a,b\right) .
\end{equation}%
\textit{If }$\alpha \geq 0$\textit{, we may allow compact subsets of }$%
[0,\infty )$\textit{.}\newline
\newline
\textbf{Remarks}\newline
(a) We remind the reader that $\mu $ is regular if $w$ is positive a.e. in $%
\left( -1,1\right) $, or more generally if $\mu \in \mathcal{M}$.\newline
(b) Our proof uses the fact that universality holds for the Jacobi weight $%
w^{\left( \alpha ,\beta \right) }$.\newline
(c) We can reformulate this in a way that allows $a,b$ to vary in a compact
subset of the complex plane. To do this one shows that $n^{-2\alpha
-2}K_{n}\left( 1-\frac{a}{2n^{2}},1-\frac{b}{2n^{2}}\right) $ is uniformly
bounded for $n\geq 1$ and $a,b$ in compact subsets of the plane. This can be
proved by bounding $\left( 1-x+\frac{1}{n^{2}}\right) ^{\alpha
+1/2}K_{n}\left( x,x\right) $ in $\left[ 1-\delta ,1\right] $ for some $%
\delta >0$, using Cauchy-Schwarz to bound $K_{n}\left( x,y\right) $, and
then using the maximum principle for subharmonic functions.

This paper is organised as follows. In the next section, we establish
asymptotics for Christoffel functions. In section 3, we localize, and in
section 4, we smoothe, and prove the theorem. In the sequel\thinspace\ $%
C,C_{1},C_{2},...$ denote constants independent of $n,x,\theta $. The same
symbol does not necessarily denote the same constant in different
occurences. We shall write $C=C\left( \alpha \right) $ or $C\neq C\left(
\alpha \right) $ to respectively denote dependence on, or independence of,
the parameter $\alpha $. Given measures $\mu ^{\ast }$, $\mu ^{\#}$, we use $%
K_{n}^{\ast },K_{n}^{\#}$ and $p_{n}^{\ast },p_{n}^{\#}$ to denote their
reproducing kernels and orthonormal polynomials. Similarly superscripts $%
\ast ,\#$ are used to distinguish their leading coefficients and Christoffel
functions, and the superscript $\left( \alpha ,\beta \right) $ denotes
quantities associated with the Jacobi weight $w^{\left( \alpha ,\beta
\right) }.$\newline
\newline
\textbf{Acknowledgement}\newline
This research was stimulated by the wonderful conference in honor of Percy
Deift's 60th birthday, held at Courant Institute in June 2006.

\section{Christoffel functions}

Recall that the $n$th Christoffel function for $\mu $ is 
\begin{equation*}
\lambda _{n}\left( x\right) =1/K_{n}\left( x,x\right) =\min_{\deg \left(
P\right) \leq n-1}\left( \int_{-1}^{1}P^{2}d\mu \right) /P^{2}\left(
x\right) .
\end{equation*}%
The methods used to prove the following result are well known, but I could
not find this theorem in the literature.\newline
\newline
\textbf{Theorem 2.1}\newline
\textit{Let }$\mu $\textit{\ be a regular measure on }$\left[ -1,1\right] .$%
\textit{\ Assume that for some }$\rho >0$, $\mu $\textit{\ is absolutely
continuous in }$J=\left[ 1-\rho ,1\right] $\textit{\ and in }$J$\textit{,} $%
w=hw^{\left( \alpha ,\beta \right) },$\textit{where }$\alpha ,\beta >-1$%
\textit{\ and }$h$\textit{\ is bounded above and below by positive constants
with} 
\begin{equation*}
\lim_{x\rightarrow 1-}h\left( x\right) =h\left( 1\right) >0.
\end{equation*}%
\textit{Let }$A>0$. \textit{Then uniformly for }$a\in \left[ 0,A\right] ,$ 
\begin{equation}
\lim_{n\rightarrow \infty }\lambda _{n}\left( 1-\frac{a}{2n^{2}}\right)
/\lambda _{n}^{\left( \alpha ,\beta \right) }\left( 1-\frac{a}{2n^{2}}%
\right) =h\left( 1\right) .
\end{equation}%
\textit{Moreover, uniformly for }$n\geq n_{0}\left( A\right) $ \textit{and} $%
a\in \left[ 0,A\right] ,$%
\begin{equation}
\lambda _{n}\left( 1-\frac{a}{2n^{2}}\right) \sim \lambda _{n}^{\left(
\alpha ,\beta \right) }\left( 1-\frac{a}{2n^{2}}\right) \sim n^{-\left(
2\alpha +2\right) }.
\end{equation}%
\textit{The constants implicit in }$\sim $\textit{\ do not depend on }$\rho $%
\textit{.} \newline
\textbf{Remark}\newline
The notation $\sim $ means that the ratio of the two Christoffel functions
is bounded above and below by positive constants independent of $n$ and $a.$
Our proof actually shows that if $\left\{ \varepsilon _{n}\right\} $ is any
sequence of positive numbers with limit $0,$ 
\begin{equation*}
\lambda _{n}\left( x\right) /\lambda _{n}^{\left( \alpha ,\beta \right)
}\left( x\right) =h\left( 1\right) +o\left( 1\right) ,
\end{equation*}%
uniformly for $x\in \left[ 1-\varepsilon _{n},1\right] $.\newline
\textbf{Proof}\newline
Let $\varepsilon >0$ and choose $\delta \in \left( 0,\rho \right) $ such
that 
\begin{equation}
\left( 1+\varepsilon \right) ^{-1}\leq \frac{h\left( x\right) }{h\left(
1\right) }\leq 1+\varepsilon ,\text{ }x\in \left[ 1-\delta ,1\right] .
\end{equation}%
Let us define a measure $\mu ^{\ast }$ with 
\begin{equation*}
\mu ^{\ast }=\mu \text{ in }[-1,1-\delta )
\end{equation*}%
and in $I=\left[ 1-\delta ,1\right] ,$ let $\mu ^{\ast }$ be absolutely
continuous, with absolutely continuous component $w^{\ast }$ satisfying 
\begin{equation}
w^{\ast }=w^{\left( \alpha ,\beta \right) }h\left( 1\right) \left(
1+\varepsilon \right) \text{ in }I.
\end{equation}%
Because of (2.3), $d\mu \leq d\mu ^{\ast }$, so that if $\lambda _{n}^{\ast }
$ is the $n$th Christoffel function for $\mu ^{\ast }$, we have for all $x$%
\begin{equation}
\lambda _{n}\left( x\right) \leq \lambda _{n}^{\ast }\left( x\right) .
\end{equation}%
We now find an upper bound for $\lambda _{n}^{\ast }\left( x\right) $\ for $%
x\in \lbrack 1-\delta /2,1]$. There exists $r\in \left( 0,1\right) $ such
that 
\begin{equation}
0\leq 1-\left( \frac{t-x}{2}\right) ^{2}\leq 1-r\text{ for }x\in \lbrack
1-\delta /2,1]\text{ and }t\in \lbrack -1,1-\delta ].
\end{equation}%
Choose $\eta \in \left( 0,\frac{1}{2}\right) $ and $\sigma >1$ so close to $1
$ that 
\begin{equation}
\sigma ^{1-\eta }<\left( 1-r\right) ^{-\eta /4}.
\end{equation}%
Let $m=m\left( n\right) =n-2\left[ \eta n/2\right] $. Fix $x\in \lbrack
1-\delta /2,1]$ and choose a polynomial $P_{m}$ of degree $\leq $ $m-1$ such
that 
\begin{equation*}
\lambda _{m}^{\left( \alpha ,\beta \right) }\left( x\right)
=\int_{-1}^{1}P_{m}^{2}w^{\left( \alpha ,\beta \right) }\text{ and }%
P_{m}^{2}\left( x\right) =1.
\end{equation*}%
\newline
Thus $P_{m}$ is the minimizing polynomial in the Christoffel function for
the Jacobi weight $w^{\left( \alpha ,\beta \right) }$ at $x$. Let 
\begin{equation*}
S_{n}\left( t\right) =P_{m}\left( t\right) \left( 1-\left( \frac{t-x}{2}%
\right) ^{2}\right) ^{\left[ \eta n/2\right] },
\end{equation*}%
a polynomial of degree $\leq $ $m-1+2\left[ \eta n/2\right] \leq n-1$ with $%
S_{n}\left( x\right) =1$. Then using (2.4) and (2.6), 
\begin{eqnarray*}
\lambda _{n}^{\ast }\left( x\right)  &\leq &\int_{-1}^{1}S_{n}^{2}d\mu
^{\ast } \\
&\leq &h\left( 1\right) (1+\varepsilon )\int_{1-\delta
}^{1}P_{m}^{2}w^{\left( \alpha ,\beta \right) }+\left\Vert P_{m}\right\Vert
_{L_{\infty }\left[ -1,1-\delta \right] }^{2}\left( 1-r\right) ^{\left[ \eta
n/2\right] }\int_{-1}^{1-\delta }d\mu ^{\ast } \\
&\leq &h\left( 1\right) (1+\varepsilon )\lambda _{m}^{\left( \alpha ,\beta
\right) }\left( x\right) +\left\Vert P_{m}\right\Vert _{L_{\infty }\left[
-1,1\right] }^{2}\left( 1-r\right) ^{\left[ \eta n/2\right]
}\int_{-1}^{1}d\mu ^{\ast }.
\end{eqnarray*}%
Now we use the key idea from \cite[Lemma 9, p. 450]{Mateetal1991}. For $%
m\geq m_{0}\left( \sigma \right) $, we have 
\begin{eqnarray*}
\left\Vert P_{m}\right\Vert _{L_{\infty }\left[ -1,1\right] }^{2} &\leq
&\sigma ^{m}\int_{-1}^{1}P_{m}^{2}w^{\left( \alpha ,\beta \right) } \\
&=&\sigma ^{m}\lambda _{m}^{\left( \alpha ,\beta \right) }\left( x\right) .
\end{eqnarray*}%
(This holds more generally for any polynomial $P$ of degree $\leq m-1$, and
is a consequence of the regularity of the measure $w^{\left( \alpha ,\beta
\right) }$. Alternatively, we could use classic bounds for the Christoffel
functions for Jacobi weight.) Then from (2.7), uniformly for $x\in \left[
1-\delta /2,1\right] ,$%
\begin{eqnarray*}
\lambda _{n}^{\ast }\left( x\right)  &\leq &h\left( 1\right) (1+\varepsilon
)\lambda _{m}^{\left( \alpha ,\beta \right) }\left( x\right) \left\{ 1+C
\left[ \sigma ^{1-\eta }\left( 1-r\right) ^{\eta /2}\right] ^{n}\right\}  \\
&\leq &h\left( 1\right) (1+\varepsilon )\lambda _{m}^{\left( \alpha ,\beta
\right) }\left( x\right) \left\{ 1+o\left( 1\right) \right\} ,
\end{eqnarray*}%
so as $\lambda _{n}\leq \lambda _{n}^{\ast },$%
\begin{eqnarray}
&&\sup_{x\in \left[ 1-\delta /2,1\right] }\lambda _{n}\left( x\right)
/\lambda _{n}^{\left( \alpha ,\beta \right) }\left( x\right)   \notag \\
&\leq &h\left( 1\right) (1+\varepsilon )\left\{ 1+o\left( 1\right) \right\}
\sup_{x\in \left[ 1-\delta /2,1\right] }\lambda _{m}^{\left( \alpha ,\beta
\right) }\left( x\right) /\lambda _{n}^{\left( \alpha ,\beta \right) }\left(
x\right) .
\end{eqnarray}%
Now for large enough $n$, and some $C$ independent of $\delta ,\eta ,m,n,$%
\begin{equation}
\sup_{x\in \left[ 1-\delta /2,1\right] }\lambda _{m}^{\left( \alpha ,\beta
\right) }\left( x\right) /\lambda _{n}^{\left( \alpha ,\beta \right) }\left(
x\right) \leq 1+C\eta .
\end{equation}%
Indeed if $\left\{ p_{k}^{\left( \alpha ,\beta \right) }\right\} $ denote
the orthonormal Jacobi polynomials for $w^{\left( \alpha ,\beta \right) }$,
they admit the bound \cite[p.170]{Nevai1979} 
\begin{equation*}
\left\vert p_{k}^{\left( \alpha ,\beta \right) }\left( x\right) \right\vert
\leq C\left( 1-x+\frac{1}{k^{2}}\right) ^{-\alpha /2-1/4},\text{ }x\in \left[
0,1\right] .
\end{equation*}%
Then 
\begin{eqnarray*}
0 &\leq &1-\frac{\lambda _{n}^{\left( \alpha ,\beta \right) }\left( x\right) 
}{\lambda _{m}^{\left( \alpha ,\beta \right) }\left( x\right) }=\lambda
_{n}^{\left( \alpha ,\beta \right) }\left( x\right)
\sum_{k=m}^{n-1}p_{k}^{2}\left( x\right)  \\
&\leq &C\lambda _{n}^{\left( \alpha ,\beta \right) }\left( x\right) \left(
n-m\right) \max_{\frac{n}{2}\leq k\leq n}\left( 1-x+\frac{1}{k^{2}}\right)
^{-\alpha -1/2} \\
&\leq &C\eta n\lambda _{n}^{\left( \alpha ,\beta \right) }\left( x\right)
\left( 1-x+\frac{1}{n^{2}}\right) ^{-\alpha -1/2} \\
&\leq &C\eta ,
\end{eqnarray*}%
by classical bounds for Christoffel functions \cite[p. 108, Lemma 5]%
{Nevai1979}. So we have (2.9). \newline
\newline
Now let $a\in \left[ 0,A\right] $. We see that for $n\geq n_{0}\left(
A\right) ,$ we have $1-a/\left( 2n^{2}\right) \in \lbrack 1-\delta /2,1]$,
and hence (2.8) gives%
\begin{equation*}
\limsup_{n\rightarrow \infty }\left( \sup_{a\in \left[ 0,A\right] }\lambda
_{n}\left( 1-\frac{a}{2n^{2}}\right) /\lambda _{n}^{\left( \alpha ,\beta
\right) }\left( 1-\frac{a}{2n^{2}}\right) \right) \leq h\left( 1\right)
(1+\varepsilon )\left( 1+C\eta \right) .
\end{equation*}%
As the left-hand side is independent of the parameters $\varepsilon ,\eta $,
we deduce that 
\begin{equation}
\limsup_{n\rightarrow \infty }\left( \sup_{a\in \left[ 0,A\right] }\lambda
_{n}\left( 1-\frac{a}{2n^{2}}\right) /\lambda _{n}^{\left( \alpha ,\beta
\right) }\left( 1-\frac{a}{2n^{2}}\right) \right) \leq h\left( 1\right) .
\end{equation}%
In a similar way, we can establish the converse bound%
\begin{equation}
\limsup_{n\rightarrow \infty }\left( \sup_{a\in \left[ 0,A\right] }\lambda
_{n}^{\left( \alpha ,\beta \right) }\left( 1-\frac{a}{2n^{2}}\right)
/\lambda _{n}\left( 1-\frac{a}{2n^{2}}\right) \right) \leq h\left( 1\right)
^{-1}.
\end{equation}%
Indeed with $m,x$ and $\eta $ as above, let us choose a polynomial $P$ of
degree $\leq m-1$ such that%
\begin{equation*}
\lambda _{m}\left( x\right) =\int_{-1}^{1}P_{m}^{2}\left( t\right) d\mu
\left( t\right) \text{ and }P_{m}^{2}\left( x\right) =1.
\end{equation*}%
Then with $S_{n}$ as above, and proceeding as above,%
\begin{equation*}
\lambda _{n}^{\left( \alpha ,\beta \right) }\left( x\right) \leq
\int_{-1}^{1}S_{n}^{2}w^{\left( \alpha ,\beta \right) }
\end{equation*}%
\begin{eqnarray*}
&\leq &\left[ h\left( 1\right) ^{-1}(1+\varepsilon )\right] \int_{1-\delta
}^{1}P_{m}^{2}d\mu +\left\Vert P_{m}\right\Vert _{L_{\infty }\left[
-1,1-\delta \right] }^{2}\left( 1-r\right) ^{\left[ \eta n/2\right]
}\int_{-1}^{1-\delta }w^{\left( \alpha ,\beta \right) } \\
&\leq &\left[ h\left( 1\right) ^{-1}(1+\varepsilon )\right] \lambda
_{m}\left( x\right) \left\{ 1+C\left[ \sigma ^{1-\eta }\left( 1-r\right)
^{\eta /2}\right] ^{n}\right\} ,
\end{eqnarray*}%
and so as above, 
\begin{eqnarray*}
&&\sup_{x\in \left[ 1-\delta /2,1\right] }\lambda _{m}^{\left( \alpha ,\beta
\right) }\left( x\right) /\lambda _{m}\left( x\right)  \\
&\leq &\left[ h\left( 1\right) ^{-1}(1+\varepsilon )(1+o\left( 1\right) )%
\right] \sup_{x\in \left[ 1-\delta /2,1\right] }\lambda _{m}^{\left( \alpha
,\beta \right) }\left( x\right) /\lambda _{n}^{\left( \alpha ,\beta \right)
}\left( x\right)  \\
&\leq &\left[ h\left( 1\right) ^{-1}(1+\varepsilon )\right] \left\{
1+o\left( 1\right) \right\} \left( 1+C\eta \right) .
\end{eqnarray*}%
Then (2.11) follows after a scale change $m\rightarrow n$ and using
monotonicity of $\lambda _{n}$ in \ $n$. Together (2.10) and (2.11) give the
result. $\blacksquare $

\section{Localization}

\textbf{Theorem 3.1}\newline
\textit{Assume that }$\mu ,\mu ^{\ast }$\textit{\ are regular measures on }$%
\left[ -1,1\right] $\textit{. Assume that }%
\begin{equation*}
d\mu =d\mu ^{\ast }=\left( hw^{\left( \alpha ,\beta \right) }\right) \left(
t\right) dt\text{ in }J=\left[ 1-\rho ,1\right] ,
\end{equation*}%
\textit{where }$h$\textit{\ satisfies the hypothesis of Theorem 2.1. Let }$%
A>0$\textit{. Then as }$n\rightarrow \infty ,$ 
\begin{equation}
\sup_{a,b\in \left[ 0,A\right] }\left\vert \left( K_{n}-K_{n}^{\ast }\right)
\left( 1-\frac{a}{2n^{2}},1-\frac{b}{2n^{2}}\right) \right\vert /n^{2\alpha
+2}=o\left( 1\right) .
\end{equation}%
\textbf{Proof\newline
}We initially assume that%
\begin{equation}
d\mu \leq d\mu ^{\ast }\text{ in }\left( -1,1\right) .
\end{equation}%
The idea is to estimate the $L_{2}$ norm of $K_{n}\left( x,t\right)
-K_{n}^{\ast }\left( x,t\right) $ over $\left[ -1,1\right] $, and then to
use Christoffel function estimates. Now 
\begin{eqnarray*}
&&\int_{-1}^{1}\left( K_{n}\left( x,t\right) -K_{n}^{\ast }\left( x,t\right)
\right) ^{2}d\mu \left( t\right)  \\
&=&\int_{-1}^{1}K_{n}^{2}\left( x,t\right) d\mu \left( t\right)
-2\int_{-1}^{1}K_{n}\left( x,t\right) K_{n}^{\ast }\left( x,t\right) d\mu
\left( t\right) +\int_{-1}^{1}K_{n}^{\ast 2}\left( x,t\right) d\mu \left(
t\right)  \\
&=&K_{n}\left( x,x\right) -2K_{n}^{\ast }\left( x,x\right)
+\int_{-1}^{1}K_{n}^{\ast 2}\left( x,t\right) d\mu \left( t\right) ,
\end{eqnarray*}%
by the reproducing kernel property. As $d\mu \leq d\mu ^{\ast },$ we also
have%
\begin{equation*}
\int_{-1}^{1}K_{n}^{\ast 2}\left( x,t\right) d\mu \left( t\right) \leq
\int_{-1}^{1}K_{n}^{\ast 2}\left( x,t\right) d\mu ^{\ast }\left( t\right)
=K_{n}^{\ast }\left( x,x\right) .
\end{equation*}%
So 
\begin{eqnarray}
&&\int_{-1}^{1}\left( K_{n}\left( x,t\right) -K_{n}^{\ast }\left( x,t\right)
\right) ^{2}d\mu \left( t\right)   \notag \\
&\leq &K_{n}\left( x,x\right) -K_{n}^{\ast }\left( x,x\right) .
\end{eqnarray}%
Next for any polynomial $P$ of degree $\leq n-1$, we have the Christoffel
function estimate 
\begin{equation}
\left\vert P\left( y\right) \right\vert \leq K_{n}\left( y,y\right)
^{1/2}\left( \int_{-1}^{1}P^{2}d\mu \right) ^{1/2}.
\end{equation}%
Applying this to $P\left( t\right) =K_{n}\left( x,t\right) -K_{n}^{\ast
}\left( x,t\right) $ and using (3.3) gives 
\begin{eqnarray*}
&&\left\vert K_{n}\left( x,y\right) -K_{n}^{\ast }\left( x,y\right)
\right\vert  \\
&\leq &K_{n}\left( y,y\right) ^{1/2}\left[ K_{n}\left( x,x\right)
-K_{n}^{\ast }\left( x,x\right) \right] ^{1/2}
\end{eqnarray*}%
so%
\begin{eqnarray*}
&&\left\vert K_{n}\left( x,y\right) -K_{n}^{\ast }\left( x,y\right)
\right\vert /K_{n}\left( x,x\right)  \\
&\leq &\left( \frac{K_{n}\left( y,y\right) }{K_{n}\left( x,x\right) }\right)
^{1/2}\left[ 1-\frac{K_{n}^{\ast }\left( x,x\right) }{K_{n}\left( x,x\right) 
}\right] ^{1/2}.
\end{eqnarray*}%
Now we set $x=1-\frac{a}{2n^{2}}$ and $y=1-\frac{b}{2n^{2}}$, where $a,b\in %
\left[ 0,A\right] $. By Theorem 2.1, uniformly for such $x,$ $\frac{%
K_{n}^{\ast }\left( x,x\right) }{K_{n}\left( x,x\right) }=1+o\left( 1\right) 
$, for they both have the same asymptotics as for the Jacobi weight.
Moreover, uniformly for $a,b\in \left[ 0,A\right] ,$ 
\begin{equation*}
K_{n}\left( 1-\frac{b}{2n^{2}},1-\frac{b}{2n^{2}}\right) \sim K_{n}\left( 1-%
\frac{a}{2n^{2}},1-\frac{a}{2n^{2}}\right) \sim n^{2\alpha +2},
\end{equation*}%
so 
\begin{equation*}
\sup_{a,b\in \left[ 0,A\right] }\left\vert \left( K_{n}-K_{n}^{\ast }\right)
\left( 1-\frac{a}{2n^{2}},1-\frac{b}{2n^{2}}\right) \right\vert /n^{2\alpha
+2}=o\left( 1\right) .
\end{equation*}%
Now we drop the extra hypothesis (3.2). Define a measure $\nu $ by $\nu =\mu
=\mu ^{\ast }$ in $J$ and 
\begin{equation*}
d\nu \left( x\right) =\max \left\{ 1,w,w^{\ast }\right\} dx+d\mu _{s}+d\mu
_{s}^{\ast },\text{ in }\left[ -1,1\right] \backslash J
\end{equation*}%
where $w,w^{\ast }$ and $\mu _{s},\mu _{s}^{\ast }$ are respectively the
absolutely continuous and singular components of $\mu ,\mu ^{\ast }$. Then $%
d\mu \leq d\nu $ and $d\mu ^{\ast }\leq d\nu $, and $\nu $ is regular as its
absolutely continuous component is positive in $\left( -1,1\right) $, and
hence lies in the even smaller class $\mathcal{M}.$ The case above shows
that the reproducing kernels for $\mu $ and $\mu ^{\ast }$ have the same
asymptotics as that for $\nu $, in the sense of (3.1), and hence the same
asymptotics as each other. $\blacksquare $

\section{Smoothing}

In this section, we approximate $\mu $ of Theorem 1.2 by a Jacobi measure $%
\mu ^{\#}$ and then prove Theorem 1.2. Our smoothing result is:\newline
\newline
\textbf{Theorem 4.1}\newline
\textit{Let }$\mu $\textit{\ be as in Theorem 1.2. Let }$\varepsilon \in
\left( 0,\frac{1}{2}\right) $\textit{\ and choose }$\delta >0$\textit{\ such
that (2.3) holds. Let} 
\begin{equation}
w^{\#}=h\left( 1\right) w^{\left( \alpha ,\beta \right) }\text{ in }\left(
-1,1\right) .
\end{equation}%
\textit{Let }$A>0$\textit{. Then there exists }$C$ and $n_{0}$\textit{\ such
that for }$n\geq n_{0},$%
\begin{equation}
\sup_{a,b\in \left[ 0,A\right] }\left\vert \left( K_{n}-K_{n}^{\#}\right)
\left( 1-\frac{a}{2n^{2}},1-\frac{b}{2n^{2}}\right) \right\vert /n^{2\alpha
+2}\leq C\varepsilon ^{1/2},
\end{equation}%
\textit{where }$C$\textit{\ is independent of }$\varepsilon ,n$.\newline
\textbf{Proof}\newline
We note that because of our localization result Theorem 3.1, we may replace $%
w$ by $w^{\ast }$, where 
\begin{equation*}
w^{\ast }=w=w^{\left( \alpha ,\beta \right) }h\text{ in }I=\left[ 1-\delta ,1%
\right] 
\end{equation*}%
and 
\begin{equation*}
w^{\ast }=w^{\left( \alpha ,\beta \right) }h\left( 1\right) \text{ in }\left[
-1,1\right] \backslash I,
\end{equation*}%
without affecting the asymptotics for $K_{n}\left( 1-\frac{a}{2n^{2}},1-%
\frac{b}{2n^{2}}\right) $. (Note that $\varepsilon $ and $\delta $ play no
role in Theorem 3.1). So in the sequel, we assume that $w=w^{\left( \alpha
,\beta \right) }h\left( 1\right) =w^{\#}$ in $\left[ -1,1\right] \backslash I
$, while keeping $w$ the same in $I$. Observe that (2.3) implies that 
\begin{equation}
\left( 1+\varepsilon \right) ^{-1}\leq \frac{w}{w^{\#}}\leq 1+\varepsilon 
\text{, in }\left[ -1,1\right] .
\end{equation}%
Then, much as in the previous section, 
\begin{eqnarray*}
&&\int_{-1}^{1}\left( K_{n}\left( x,t\right) -K_{n}^{\#}\left( x,t\right)
\right) ^{2}w^{\#}\left( t\right) dt \\
&=&\int_{-1}^{1}K_{n}^{2}\left( x,t\right) w^{\#}\left( t\right)
dt-2\int_{-1}^{1}K_{n}\left( x,t\right) K_{n}^{\#}\left( x,t\right)
w^{\#}\left( t\right) dt+\int_{-1}^{1}K_{n}^{\#2}\left( x,t\right)
w^{\#}\left( t\right) dt \\
&=&\int_{-1}^{1}K_{n}^{2}\left( x,t\right) w\left( t\right)
dt+\int_{I}K_{n}^{2}\left( x,t\right) \left( w^{\#}-w\right) \left( t\right)
dt-2K_{n}\left( x,x\right) +K_{n}^{\#}\left( x,x\right)  \\
&=&K_{n}^{\#}\left( x,x\right) -K_{n}\left( x,x\right)
+\int_{I}K_{n}^{2}\left( x,t\right) \left( w^{\#}-w\right) \left( t\right)
dt,
\end{eqnarray*}%
recall that $w=w^{\#}$ in $\left[ -1,1\right] \backslash I$. By (4.3), 
\begin{equation*}
\int_{I}K_{n}^{2}\left( x,t\right) \left( w^{\#}-w\right) \left( t\right)
dt\leq \varepsilon \int_{J}K_{n}^{2}\left( x,t\right) w\left( t\right)
dt\leq \varepsilon K_{n}\left( x,x\right) .
\end{equation*}%
So%
\begin{equation}
\int_{-1}^{1}\left( K_{n}\left( x,t\right) -K_{n}^{\#}\left( x,t\right)
\right) ^{2}w^{\#}\left( t\right) dt\leq K_{n}^{\#}\left( x,x\right) -\left(
1-\varepsilon \right) K_{n}\left( x,x\right) .
\end{equation}%
Applying an obvious analogue of (3.4) to $P\left( t\right) =K_{n}\left(
x,t\right) -K_{n}^{\#}\left( x,t\right) $ and using (4.4) gives for $y\in %
\left[ -1,1\right] ,$ 
\begin{eqnarray*}
&&\left\vert K_{n}\left( x,y\right) -K_{n}^{\#}\left( x,y\right) \right\vert 
\\
&\leq &K_{n}^{\#}\left( y,y\right) ^{1/2}\left[ K_{n}^{\#}\left( x,x\right)
-\left( 1-\varepsilon \right) K_{n}\left( x,x\right) \right] ^{1/2}
\end{eqnarray*}%
so%
\begin{eqnarray*}
&&\left\vert K_{n}\left( x,y\right) -K_{n}^{\#}\left( x,y\right) \right\vert
/K_{n}^{\#}\left( x,x\right)  \\
&\leq &\left( \frac{K_{n}^{\#}\left( y,y\right) }{K_{n}^{\#}\left(
x,x\right) }\right) ^{1/2}\left[ 1-\left( 1-\varepsilon \right) \frac{%
K_{n}\left( x,x\right) }{K_{n}^{\#}\left( x,x\right) }\right] ^{1/2}.
\end{eqnarray*}%
In view of (4.3), we also have 
\begin{equation*}
\frac{K_{n}\left( x,x\right) }{K_{n}^{\#}\left( x,x\right) }=\frac{\lambda
_{n}^{\#}\left( x\right) }{\lambda _{n}\left( x\right) }\geq \frac{1}{%
1+\varepsilon },
\end{equation*}%
so for all $y\in \left[ -1,1\right] ,$%
\begin{eqnarray*}
&&\left\vert K_{n}\left( x,y\right) -K_{n}^{\#}\left( x,y\right) \right\vert
/K_{n}^{\#}\left( x,x\right)  \\
&\leq &\left( \frac{K_{n}^{\#}\left( y,y\right) }{K_{n}^{\#}\left(
x,x\right) }\right) ^{1/2}\left[ 1-\frac{1-\varepsilon }{1+\varepsilon }%
\right] ^{1/2} \\
&\leq &\sqrt{2\varepsilon }\left( \frac{K_{n}^{\#}\left( y,y\right) }{%
K_{n}^{\#}\left( x,x\right) }\right) ^{1/2}.
\end{eqnarray*}%
Now we set $x=1-\frac{a}{2n^{2}}$ and $y=1-\frac{b}{2n^{2}}$, where $a,b\in %
\left[ 0,A\right] $. By Theorem 2.1, uniformly for $a,b\in \left[ 0,A\right]
,$ 
\begin{equation*}
K_{n}^{\#}\left( 1-\frac{b}{2n^{2}},1-\frac{b}{2n^{2}}\right) \sim
K_{n}^{\#}\left( 1-\frac{a}{2n^{2}},1-\frac{a}{2n^{2}}\right) \sim
n^{2\alpha +2},
\end{equation*}%
and also the constants implicit in $\sim $ are independent of $\varepsilon $
(this is crucial!). Thus for some $C$ and $n_{0}$ depending only on $A$, we
have for $n\geq n_{0},$%
\begin{equation*}
\sup_{a,b\in \left[ 0,A\right] }\left\vert \left( K_{n}-K_{n}^{\#}\right)
\left( 1-\frac{a}{2n^{2}},1-\frac{b}{2n^{2}}\right) \right\vert /n^{2\alpha
+2}\leq C\sqrt{\varepsilon }.
\end{equation*}%
$\blacksquare $\newline
\newline
\textbf{Proof of Theorem 1.2}\newline
Let $\varepsilon _{1}>0$. We can choose $\varepsilon >0$ so small that the
right-hand side of (4.2) is less than $\varepsilon _{1}$. (Recall that $C$
there is independent of $\varepsilon $). Hence for $n\geq n_{0}\left(
A,\varepsilon _{1}\right) ,$%
\begin{equation*}
\sup_{a,b\in \left[ 0,A\right] }\left\vert \left( K_{n}-K_{n}^{\#}\right)
\left( 1-\frac{a}{2n^{2}},1-\frac{b}{2n^{2}}\right) \right\vert /n^{2\alpha
+2}\leq \varepsilon _{1}.
\end{equation*}%
It follows that 
\begin{equation*}
\lim_{n\rightarrow \infty }\left( \sup_{a,b\in \left[ 0,A\right] }\left\vert
\left( K_{n}-K_{n}^{\#}\right) \left( 1-\frac{a}{2n^{2}},1-\frac{b}{2n^{2}}%
\right) \right\vert /n^{2\alpha +2}\right) =0.
\end{equation*}%
Next, uniformly for $a\in \left[ A_{1},A_{2}\right] $, where $%
0<A_{1}<A_{2}<\infty $, we see that 
\begin{equation*}
w\left( 1-\frac{a}{2n^{2}}\right) =\left( 1+o\left( 1\right) \right) h\left(
1\right) 2^{\beta }\left( \frac{a}{2n^{2}}\right) ^{\alpha }=w^{\#}\left( 1-%
\frac{a}{2n^{2}}\right) \left( 1+o\left( 1\right) \right) ,
\end{equation*}%
with a similar relation when we replace $a$ by $b$. Hence uniformly for $%
a,b\in \left[ A_{1},A_{2}\right] $, 
\begin{eqnarray*}
&&\frac{1}{2n^{2}}\tilde{K}_{n}\left( 1-\frac{a}{2n^{2}},1-\frac{b}{2n^{2}}%
\right)  \\
&=&\frac{1}{2n^{2}}\tilde{K}_{n}^{\#}\left( 1-\frac{a}{2n^{2}},1-\frac{b}{%
2n^{2}}\right) \left( 1+o\left( 1\right) \right) +o\left( 1\right)  \\
&=&\mathbb{J}_{\alpha }\left( a,b\right) +o\left( 1\right) ,
\end{eqnarray*}%
by the universality limit (1.2) for the scaled Jacobi weight $w^{\#}=h\left(
1\right) w^{\left( \alpha ,\beta \right) }$. For this, see for example \cite%
{KuijlaarsVanlessen2002}. When $\alpha \geq 0$, we can allow instead $a\in %
\left[ 0,A\right] $. $\blacksquare $\newline

\end{document}